\documentclass[a4paper,10pt]{article} 
\usepackage{latexsym,amsfonts,amsmath}
\usepackage[english,francais]{babel}

\usepackage{rotating}
\usepackage{graphics,epsf,psfrag,eepic,epic} 

\author{Charles
  Cochet\footnote{
    {\sf cochet@math.jussieu.fr} or 
    http://www.math.jussieu.fr/$\sim$cochet}
    } 
\title{Multiplicities and tensor product coefficients for $A_r$} 

\newtheorem{theo}{Theorem}[section]

\newtheorem{lemm}[theo]{Lemma}

\newtheorem{rema}[theo]{Remark}

\newcommand{\Zset}{\mathbb{Z}}

\newcommand{\N}{\mathbb{N}}
\newcommand{\Z}{\mathbb{Z}}

\newcommand{\R}{\mathbb{R}}
\newcommand{\C}{\mathbb{C}}

\def\bnum{\begin{enumerate}}
\def\enum{\end{enumerate}}
\def\la{\lambda}

\def\LiE.{{\sf L\kern-.25em\raise0.59ex\hbox{\i}\kern-0.03em E}}
\def\sec{\operatorname{s}}
\def\eps{\varepsilon}

\def\clm{c_\la^\mu}
\def\clmn{c_{\la,\mu}^\nu}

\def\res{\operatorname{Res}}
\def\ires{\operatorname{IRes}}

\def\rmo{{\rm (}} 
\def\rmf{{\rm )}} 
\def\defa{\operatorname{def}(a)}
\def\spa{\operatorname{Sp}(a)}

\begin{document} 
\selectlanguage{english}
\maketitle

\abstract{
We apply some recent developments of
Baldoni-DeLoera-Vergne~\cite{BdLV} on vector partition functions, 
to Kostant and Steinberg formulas, in the case of $A_r$. 
We therefore get a fast {\sc Maple} program that computes for $A_r$: 
the multiplicity $\clm$ of the weight $\mu$ in the representation
$V(\la)$ of highest weight $\lambda$; 
the multiplicity $\clmn$ of the representation $V(\nu)$ in
$V(\la)\otimes V(\mu)$. 
The computation also gives the locally polynomial functions 
$\clm$ and $\clmn$. 
}

\section{Introduction}

In this short note, we are interested in the two following problems in
the case of $A_r$: 

\bnum 
\item[{\bf Mult:}] Computation of the multiplicity $\clm$ of the
  weight $\mu$ in the representation $V(\la)$ of highest weight
  $\la$. 
\item[{\bf Tens:}] Computation of the multiplicity $\clmn$ of the 
  representation $V(\nu)$ in the tensor product of representations of 
  highest weights $\la$ and $\mu$. 
\enum 

The approach to these problems is through vector partition 
functions, namely number of integral points in lattice polytopes. 
More precisely, let $\Phi$ be a $n\times N$ integral matrix with 
column vectors $\Phi_1$, \ldots, $\Phi_N$. 
Fix a $n$-dimensional vector $a$. 
The rational convex polytope associated to $\Phi$ and $a$ is 
$$P(\Phi,a)=\left\{x\in\R^N\,;\,
  \sum_{i=1}^N x_i\Phi_i=a,\,x_i\geq 0\right\}.$$ 
We assume that $a$ is in the cone $C(\Phi)$ spanned by 
non-negative linear combinations of the vectors $\Phi_i$. 
We also assume that $\ker(\Phi)\cap\R_+^N=\{0\}$, so that the cone 
$C(\Phi)$ is acute. 
The {\it vector partition function} is then by definition 
$$k(\Phi,a)=\left|P(\Phi,a)\cap\N^N\right|,$$ 
that is the number of non-negative integral solutions 
$(x_1,\ldots,x_N)$ of the equation $\sum_{i=1}^N x_i\Phi_i=a$. 
If $\Phi$ is the matrix of positive roots of $A_r$, 
then $k(\Phi,a)$ is denoted by $k(A_r^+,a)$ and it is called 
{\it Kostant partition function}. 
Note that $\Phi$ is the $(r+1)\times(r(r+1)/2)$ matrix with columns
$e_i-e_j$ ($1\leq i<j\leq r+1$), where $e_i$ is the canonical basis
of $\R^{r+1}$. 

\bigskip 

Let $\Sigma_{r+1}$ be the set of permutations of $(r+1)$ elements. 
This is the Weyl group of $A_r$. Kostant multiplicity formula asserts
that 
\begin{equation} \label{eq:kostant} 
\clm = \sum_w (-1)^{\eps(w)}
k\Big(A_r^+,w(\la+\rho)-(\mu+\rho)\Big), 
\end{equation} 
where $\rho$ is half the sum of positive roots. 
Here, the sum is over the elements $w\in\Sigma_{r+1}$ such that 
$w(\la+\rho)-(\mu+\rho)$ is in the cone generated by non-negative
combinations of positive roots. Moreover $\eps(w)$ is the signature 
of $w$. 

Steinberg formula asserts that 
\begin{equation} \label{eq:steinberg}
\clmn = \sum_{w,w'} (-1)^{\eps(w)\eps(w')} 
 k\Big(A_r^+,w(\la+\rho)+w'(\mu+\rho)-(\nu+2\rho)\Big). 
\end{equation} 
Here, the sum is over couples
$(w,w')\in\Sigma_{r+1}\times\Sigma_{r+1}$ such that 
$w(\la+\rho)+w'(\mu+\rho)-(\nu+2\rho)$ is in the cone $C(A_r^+)$. 

\bigskip 

We use results of Baldoni-DeLoera-Vergne~\cite{BdLV} and
Baldoni-Vergne~\cite{BalVer01} on vector partition functions to 
obtain an efficient {\sc Maple} program. 
Vector partition function is computed via inverse Laplace formula,
involving iterated residues of rational functions. 

Recall that \LiE. program (see~\cite{MvL94}) uses 
Freudenthal and Klymik formulas. 
The program \LiE. is designed to work for any root system, 
while our program is designed specially for large parameters in 
$A_r$. 

\newpage
\section{Baldoni-DeLoera-Vergne formula} 
\label{sec:BdLV} 

Consider a $r+1$ real dimensional vector space. Let $A_r^+$ (the
positive root system of $A_r$) be defined by 
$$A_r^+=\{(e_i-e_j)\,;\,1\leq i<j\leq(r+1)\}.$$
Let $E_r$ be the vector space spanned by the elements
$(e_i-e_j)$. Then 
$$E_r=\{a\in\R^{r+1}\,;\, 
a=a_1 e_1+\cdots+a_r e_r+a_{r+1} e_{r+1}\ \ \textrm{ with }\ \ 
a_1+a_2+\cdots+a_r+a_{r+1}=0\}.$$
 The vector space $E_r$ is of dimension $r$ and the map
\begin{equation}\label{eq:iden}
f: \ \R^r \longrightarrow E_r
\end{equation}
defined by
$$a=(a_1,a_2,\ldots, a_r) \longmapsto a_1 e_1+\cdots+a_{r}
e_r-(a_1+\cdots + a_r)e_{r+1}$$ 
explicitly provides an isomorphism of $E_r$ with the Euclidean space
$\R^r.$ 
The hyperplane arrangement (setting $z_{r+1}=0$) generated by
$A_r^+$ is given by the following set of hyperplanes in $\C^r$: 
$$\{z_i\,;\,1\leq i\leq r\}\cup\{(z_i-z_j)\,;\,1\leq i<j\leq r\}.$$

Let $R_{A_r}$ be the set of rational functions 
$f(z_1,z_2,\ldots,z_r)$ on $\C^r$, with poles on the hyperplanes
$z_i=z_j$ or $z_i=0$. 
For a permutation $w\in\Sigma_r$ define the linear form on $R_{A_r}$
\begin{eqnarray*} 
\ires_{z=0}^{w}f & = & 
 \res_{z_{w(1)}=0}\res_{z_{w(2)}=0}\cdots
 \res_{z_{w(r)}=0}f(z_1,z_2,\ldots, z_r)\\ 
 & = & \res_{z_1=0}\res_{z_2=0}\cdots
 \res_{z_r=0}f(z_{w^{-1}(1)},z_{w^{-1}(2)},\ldots,
  z_{w^{-1}(r)}). 
\end{eqnarray*} 
In particular for $w=\operatorname{id}$ the linear form $f\mapsto
\ires_{z=0}f$ defined by 
$$\ires_{z=0}f=\res_{z_1=0}\res_{z_2=0}\cdots
 \res_{z_r=0}f(z_1,z_2,\ldots,z_r)$$ 
is called the {\em iterated residue.}

\bigskip 

Let $C(A_r^+)\subset E_r$ be the cone generated by positive roots. 
A subset $\sigma$ of $A_r^+$ is called a {\em basic} subset if 
$\{\sigma\}$ form a vector space basis of $E_r$. 
The {\em chamber complex} is the polyhedral subdivision of 
the cone $C(A_r^+)$ which is defined as the common refinement of 
the simplicial cones $C(\sigma)$ running over all possible basic 
subsets of $A_r^+.$ The pieces of this subdivision are called 
{\em chambers}. 
See~\cite{BdLV} and~\cite{dLStu} 
for the computation of chambers for $A_r^+$. 

A {\em wall} is a hyperplane in $E_r$ spanned by $(r-1)$ vectors of
$A_r^+$. A vector $v\in C(A_r^+)$ is {\em regular} if it is not on a 
wall. This means that for every strict subset 
$I\subset\{1,\ldots,r+1\}$ we have $\sum_{i\in I}v_i\neq 0$. 

\bigskip 

Let $\spa$ be the set of permutations $w\in\Sigma_r$ such that: 
$$\left \{ \begin{array}{rlrlrl}
\mbox{if} & a_{w(1)}\geq 0 & \mbox{then} & w(1)<w(2) 
& \mbox{else} & w(1)>w(2)\\
\mbox{if} & a_{w(1)}+a_{w(2)}\geq 0 & \mbox{then} & w(2)<w(3) 
& \mbox{else} & w(2)>w(3)\\
& \cdots &&&& \\
\mbox{if} & a_{w(1)}+\cdots+a_{w(i)}\geq 0 & \mbox{then} 
& w(i)<w(i+1) & \mbox{else} & w(i)>w(i+1)\\
& \cdots &&&&\\
\mbox{if} & a_{w(1)}+\cdots+a_{w(r-1)}\geq 0 & \mbox{then}
& w(r-1)<w(r) & \mbox{else} & w(r-1)>w(r)\\
\end{array} \right\}$$ 
An element of $\spa$ will be called a {\em special permutation} for
$a$. Remark that if $a_i\geq 0$ for all $i\leq r$, then 
$\spa=\{\operatorname{id}\}$. 
Also remark that $\spa$ is a subset of the subgroup $\Sigma_r$ of the
Weyl group $\Sigma_{r+1}$ of $A_r$. 

\bigskip 

Given $a\in C(A^+_r)\cap \Z^{r+1},$ define 
$\defa=a+\eps(\sum_{i=1}^{r}e_i-re_{r+1})$ 
with $\eps=\frac{1}{2r}$. 

\begin{lemm} \label{lemm:deformation} 
The deformed vector $\defa$ verifies: 
\begin{itemize}
\item $\defa$ is regular. 
\item $a\in C(A^+_r)$ if and only if $\defa\in C(A^+_r)$. 
\end{itemize}
\end{lemm}

{\bf Proof:} 

\begin{itemize} 
\item Let $I_a$ be a strict subset of $\{1,\ldots,r+1\}$ such that 
  $\sum_{i\in I_a}\defa_i=0$. 

  First, assume that $r+1\notin I_a$. We can re-index $a$ in order to
  get $I_a=\{1,\ldots,k\}$ with $k\leq r$. 
  Thus $(a_1+\eps)+\cdots+(a_k+\eps)=0$ means that the integer
  $a_1+\cdots+a_k$ equals $-\frac{k}{2r}$. But $0<k\leq r$ implies
  $0<\frac{k}{2r}<1$, contradiction. 

  Now, assume that $r+1\in I_a$. We can also assume that 
  $I_a=\{1,\ldots,k,r+1\}$ with $k\leq r$. 
  By definition $(a_1+\eps)+\cdots+(a_k+\eps)+(a_{r+1}-r\eps)$ equals
  $0$, therefore the integer $a_1+\cdots+a_k+a_{r+1}$ is equal to
  $\frac{r-k}{2r}$. But $k\leq r$ leads to 
  $-\frac 12 < \frac{r-k}{2r} < \frac 12$, hence $k=r$. 
  Consequently $I_a=\{1,\ldots,r+1\}$, contradiction. 
\item Note that the coordinates $a_i$ of $a$ are integers. 
  Now the integer $a_1+\cdots+a_i$ is non-negative if and only if
  $a_1+\cdots+a_i+\frac{1}{2r}$ is non-negative, because
  $0<\frac{1}{2r}<1$. Hence $a\in C(A_r^+)$ is equivalent to 
  $\defa\in C(A_r^+)$. 
\end{itemize} 

\bigskip 

Now we can state the formula that was implemented: 

\begin{theo}[Baldoni-DeLoera-Vergne~\cite{BdLV}] 
\label{theo:BdLV}
For $a\in C(A^+_r)\cap\Z^{r+1},$ the Kostant partition function
is given by:
$$k(A^+_r,a)=\sum_{w\in Sp(a')}(-1)^{n(w)}\ires_{z=0}^w
\left(\frac{(1+z_1)^{a_1+r-1}(1+z_2)^{a_2+r-2}\cdots (1+z_r)^{a_r}}
     {z_1\cdots z_r\prod_{1\leq i<j\leq r}(z_i-z_{j})}\right)$$ 
where 
$$a'=\left\{\begin{array}{cc}
a & \textrm{if a is regular,}\\
\defa & \textrm{otherwise.}
\end{array}\right.$$ 
In particular, if $a_i\geq 0$ for $1\leq i\leq r$, we have
$$k(A^+_r,a)=\res_{z_1=0}\res_{z_2=0}\cdots\res_{z_r=0}
\left(
\frac{(1+z_1)^{a_1+r-1}(1+z_2)^{a_2+r-2}\cdots(1+z_r)^{a_r}}
     {z_1\cdots z_r\prod_{1\leq i<j\leq r}(z_i-z_{j})}\right).$$
\end{theo}

\section{{\it Deus ex machina}}

This section features a brief description of the algorithms 
that were implemented with the software {\sc Maple}. 
This program is available at 
{\sf http://www.math.jussieu.fr/$\sim$cochet}

\subsection{How to use the program}

The initial data are only vectors: 
two for computing the multiplicity $\clm$, 
three for computing the tensor product coefficient $\clmn$. 

\bigskip 

Our program works with weights represented in the canonical basis of 
$\R^{r+1}$, and not fundamental weights basis of $A_r$ like \LiE.. 
The translation between these two approaches is performed via the procedures 
{\sf fundamental} and {\sf fundamental\_inverse}. 
For example {\sf fundamental([2,1,-3])} returns {\sf [1,4]}. 

Therefore computing the multiplicity $\clm$ is done by typing in 
{\sf multiplicity(lambda,mu)} where $\la$ and $\mu$ are lists of $r+1$ 
rationals such that $\sum_{i=1}^{r+1}\la_i=\sum_{i=1}^{r+1}\mu_i$ 
and $\la_i-\la_{i+1}\in\N$, $\mu_i-\mu_{i+1}\in\Zset$. 

For computing the tensor product coefficient $\clmn$, the syntax is 
{\sf tensor\_product(lambda,mu,nu)} where 
$\la$, $\mu$ and $\nu$ are lists of $r+1$ rationals such that 
$\sum_{i=1}^{r+1}(\la_i+\mu_i)=\sum_{i=1}^{r+1}\nu_i$ and 
$\la_i-\la_{i+1}\in\N$, $\mu_i-\mu_{i+1}\in\N$,
$\nu_i-\nu_{i+1}\in\N$. 



\bigskip 

In the examples, we use the vector 
$\theta=r e_1+(r-1)e_2+\cdots+1 e_{r-1}-\frac{r(r+1)}{2}e_{r+1}$. 
Its decomposition in the fundamental weights basis is the 
$r$-dimensional vector $(1,\ldots,1,1+r(r+1)/2)$. 


\subsection{Implementation} 

The elements we need to compute are: 
\bnum 
\item The vector $a'=\defa$ obtained by deforming the initial
  parameter $a$. 
\item The residues that appear in theorem~\ref{theo:BdLV}. 
\item The two sets of permutations that appear in Kostant and
  Steinberg formulas (see~\eqref{eq:kostant}
  and~\eqref{eq:steinberg}). 
\item The set of special permutations $\operatorname{Sp}(\defa)$. 
\enum 

Because of lemma~\ref{lemm:deformation}, we may use $\defa$ instead of
$a$ and we do this to simplify the procedures. 
We compute the vector $\defa$ via the straightforward {\sc Maple}
procedure {\sf defvector}. This takes care of the first part. 

\bigskip 

Computation of residues is done iteratively. The function $F$ which
residue we need to compute is a product of a certain number of
functions. This allows to take the residues by introducing little by
little the part of the function $F$ containing the needed
variable. See a detailed explanation of this procedure
in~\cite{BdLV}.  

\bigskip 

Let $u$, $v\in E^r$. A {\it valid permutation} for $u$ and $v$ is a
permutation $w\in\Sigma_{r+1}$ such that $w(u)-v\in C(A_r^+)$. 
We denote by $V(u,v)$ the set of valid permutations for $u$, $v$. 
Hence, Kostant formula for $A_r$ rewrites as 
$$\clm=\sum_{w\in V(\la+\rho,\mu+\rho)} 
  (-1)^{\eps(w)}\ k(A_r^+,w(\la+\rho)-(\mu+\rho)).$$ 

Given a set of chambers $\{C_w\}_{w\in V(\la,\mu)}$ of $C(A_r^+)$, it
follows from~\cite{SzeVer02} that $\clm$ is polynomial when 
$w\la-\mu\in\overline{C_w}$, for $w\in V(\la,\mu)$. 
In particular, the function $N\mapsto c_{N\la,N\mu}$ is a polynomial
in $N$ of degree less of equal to $\frac{r(r-1)}{2}$. 

\bigskip 

Let us explain our implementation with the symbolic langage 
{\sc Maple} of the procedure {\sf valid\_permutations} designed to
find the set $V(u,v)$. 
The method is quite simple: we build the permutations iteratively. 
This allows us not examining all permutations and saving much time. 
Recall that we have to find all permutations $w$'s such that 
$u_{w(1)}\geq v_1$, $u_{w(1)}+u_{w(2)}\geq v_1+v_2$, etc. 
For any sequence $x$ of indices, we denote by $u_x$ the sum
$\sum_{i\in x}u_i$. 

\bnum 
\item[Step $1$.] Let $X$ be the set of all indices $i$ such that
 $u_i\geq v_1$. 
\item[Step $2$.] For each $x\in X$, we find all indices 
 $i_x$ such that $u_x+u_{i_x}\geq v_x+v_{i_x}$. 
 Let $X_{new}$ be the set of such $[x,i_x]$, for all $x$ and $i_x$. 
 Then $X\leftarrow X_{new}$. 
\enum 
We repeat $r$ times step $2$, and obtain the list $X$ of 
$(r+1)$-uples representing permutations of $1$, \ldots, $r+1$. 
The second step is treated in the procedure 
{\sf next\_index\_valid\_permutations}. 
The procedure {\sf valid\_permutations} contains first step and a 
{\sf for\ldots do} loop executing $r$ times step $2$. 

\begin{rema} 
We reduce computing time by using the following three tricks. 
\bnum 
\item We compute once and for all the vector 
 $v'=[v_1,v_1+v_2,\ldots,v_1+\cdots+v_{r+1}]$. 
\item We build at the same time of 
 $X=[[i_1,\ldots,i_p],\ldots,\mbox{\rmo other sets of indices\rmf}]$ 
 the set $SX$ of partial sums associated to each
 $[i_1\ldots,i_p]$. More precisely 
 $SX=[u_{i_1}+\cdots+u_{i_p},\ldots,
   \mbox{\rmo other partial sums\rmf}]$. 
\item We use tables instead of lists. 
\enum 
\end{rema} 

Now let us examine the couples of permutations involved in Steinberg
formula. Let $u_1$, $u_2$, $v\in E^r$. 
A {\it valid couple of permutations} for $u_1$, $u_2$ and $v$ is a
couple $(w_1,w_2)\in\Sigma_{r+1}\times\Sigma_{r+1}$ such that 
$w_1(u_1)+w_2(u_2)-v\in C(A_r^+)$.  
We denote by $V(u_1,u_2,v)$ the set of valid couples of permutations
for $u_1$, $u_2$ and $v$. Hence Steinberg formula rewrites as 
$$\clmn = \sum_{(w,w')\in V(\la+\rho,\mu+\rho,\nu+2\rho)} (-1)^{\eps(w)\eps(w')} 
  \ k(A_r^+,(w(\la+\rho)+w'(\mu+\rho)-(\nu+2\rho)).$$ 
The procedure computing valid couples of permutations is similar
to the former. 

\bigskip 

To compute the subset $\spa$ of $\Sigma_r$, we use the procedure 
{\sf special\_permutations}. This procedure is very similar to the
previous one. We stress that the {\sc Maple} function 
{\sf combinat[permute]} is impractical and does not go very far
because of memory limitations. 

\section{Test of the program}

Let $\theta$ be the $r$-dimensional vector $(1,\ldots,1,1+r(r+1)/2)$
(fundamental weights decomposition in $A_r$). 
It translates as $(r,r-1,\ldots,1,-r(r+1)/2)$ in the canonical 
basis of $\R^{r+1}$. 
We used this vector to check the well-known fact that the multiplicity
of the weight $0$ in the representation of $A_r$ of highest weight 
$N\theta$ is given by the dimension of the representation of $A_{r-1}$ 
of highest weight $N\rho$, which is $(N+1)^{r(r-1)/2}$. 

\bigskip 

In this test, we compute for various $A_r$ ($r=1$, \ldots, $8$): 
\begin{itemize} 
\item $\clm$ with $\la=N\theta$ and either $\mu=0$ (worst case),
  or $\mu=[9N/10]\theta$ (intermediate case). 
\item $\clmn$ with $\la=\mu=N\theta$ and either $\nu=0$ (worst case),
  or $\nu=2[9N/10]\theta$ (intermediate case). 
\end{itemize} 

Tests were made with bi-processor PIII $1,13$GHz. 
The notation ''$-$'' in an array means that we did not try the
computation (and not that computation failed). 

\bigskip 

Recall (see for example~\cite{Bar97} and~\cite{BarPom99}) that
counting integral points in a lattice polytope is 
polynomial in the size of input if dimension is fixed, 
and NP-hard if dimension is not fixed. 
The figures~\ref{poly1} and~\ref{poly2} emphasizes this result. 

\bigskip 


In figure~\ref{poly1}, the letter I stands for intermediate case 
($\mu=[9N/10]\theta$ for $\clm$ and $\nu=2[9N/10]\theta$ for $\clmn$), 
while W stands for worst case 
($\mu=0$ for $\clm$ and $\nu=0$ for $\clmn$).

\begin{figure}[htb] 
\begin{center} 
$$\begin{array}{|c|c|c|c|c|c|c|c|c|c|} \hline 
& N=10^1 & N=10^2 & N=10^3 & N=10^4 & N=10^5 & 
N=10^6 & N=10^7 & N=10^8 & N=10^9\\ \hline 
\textrm{$\clm$, I, $A_7$} 
 & 12.5\sec & 16.0\sec & 17.5\sec & 18.7\sec & 17.63\sec
 & 19.4\sec & 20.3\sec & 21.3\sec & 22.5 \\ \hline 
\textrm{$\clm$, W, $A_7$}
 & 204.9\sec & 221.0\sec & 235.6\sec & 251.5\sec & 259.1\sec 
 & 261.5\sec & 283.8\sec & 297.0\sec & 297.2\sec\\ \hline 
\textrm{$\clmn$, I, $A_6$}
 & 40.5\sec & 47.0\sec & 50.3\sec & 52.6\sec & 53.8\sec 
 & 57.0\sec & 58.4\sec & 59.9\sec & 62.3\sec\\ \hline 
 \textrm{$\clmn$, W, $A_4$}
 & 13.5\sec & 13.7\sec & 13.8\sec & 14.0\sec & 14.1\sec 
 & 14.4\sec & 15.1\sec & 15.2\sec & 15.5\sec\\ \hline 
\end{array}$$
\end{center} 
\caption{Time of computation, when size of input grows} 
\label{poly1}
\end{figure}


\newpage 

The computation can also be done with parameters, giving
$(N+1)^{r(r-1)/2}$ as expected. 

\begin{figure}[htb] 
\begin{center} 
$$\begin{array}{|c||c|c||c|c|} \hline 
\textrm{Algebra} & \textrm{Time} & \textrm{Multiplicity } c_\theta^0 & 
\textrm{Time} & \textrm{Polynomial $N\mapsto c_{N\theta}^0$}\\
\hline\hline  
A_2 & <0.1\sec & 2=2^1 & <0.1\sec & (N+1)^1\\ \hline
A_3 & <0.1\sec & 8=2^3 & <0.1\sec & (N+1)^3\\ \hline
A_4 & <0.1\sec & 64=2^6 & <0.1\sec & (N+1)^6\\ \hline
A_5 & 0.4\sec & 1024=2^{10} & 1.4\sec & (N+1)^{10}\\ \hline
A_6 & 7.6\sec & 32768=2^{15} & 36.2\sec & (N+1)^{15}\\ \hline
A_7 & 169.3\sec & 2097152=2^{21} & 2091\sec & (N+1)^{21}\\ \hline
A_8 & 9401\sec & 268435456=2^{28} & - & -\\ \hline
\end{array}$$
\end{center} 
\caption{Multiplicity of $0$ in $V(N\theta)$ when rank increases} 
\label{poly2}
\end{figure}


\bigskip 

{\bf Acknowledgements:}
Marc A. A. van Leeuwen explained me his software \LiE. 
and its internal mechanisms. 
Moreover, his competence in {\sc Maple} has deeply influenced my
method of programing.


\end{document}